\documentclass[12pt]{amsart}
\usepackage{amscd,amssymb,righttag}
\newtheorem{thm}{Theorem}[section]
\newtheorem{corol}[thm]{Corollary} 
\newtheorem{lemma}[thm]{Lemma}
\newtheorem{prop}[thm]{Proposition}
\theoremstyle{definition}
\newtheorem{defin}[thm]{Definition}

\theoremstyle{remark}
\newtheorem{remark}[thm]{Remark}
\newtheorem{example}[thm]{Example}

\numberwithin{equation}{section}

\newcommand{\abs}[1]{\lvert#1\rvert}
\def\norm#1{\left\Vert#1\right\Vert}
\newcommand{\tri}{\hfill$\blacktriangle$} 
\def\R {{\mathbb R}}
\def\C {{\mathbb C}}
\def\N{{\mathbb N}}
\def\e{{\epsilon}}
\def\Z {{\mathbb Z}}
\def\s{{\mathbb S}}
\def\H{{\mathcal H}}

\begin{document}

\title[Dynamics of the unit sphere]{Amenable representations 
and dynamics of the unit sphere in an infinite-dimensional Hilbert space}

\author[V.G. Pestov]{Vladimir G. Pestov}
\address{School of Mathematical and Computing Sciences,
Victoria University of Wellington, P.O. Box 600, Wellington,
New Zealand}
\email{vova@mcs.vuw.ac.nz}
\urladdr{http://www.vuw.ac.nz/$^\sim$vova}
\thanks{This research was partially supported by the Marsden Fund
grant \#VUW703 of the Royal Society of New Zealand.}
\thanks{Preprinted as: Reseach Report 99-10, School of Mathematical and
Computing Sciences, Victoria University of Wellington, March 1999.}

\date{March 15, 1999}



\begin{abstract} We establish a close link between the amenability
property of a unitary representation $\pi$ of a group $G$ 
(in the sense of Bekka)
and the concentration property (in the sense of
V. Milman) of the corresponding dynamical system $(\s_\pi,G)$, where 
$\s_\H$ is the unit sphere the Hilbert space of representation.
We prove that $\pi$ is amenable if and only if either
$\pi$ contains a finite-dimensional subrepresentation or
the maximal uniform compactification of $\s_\pi$
has a $G$-fixed point. Equivalently, the latter means that the $G$-space
$(\s_\pi,G)$ has the concentration property:
every finite cover of the sphere $\s_\pi$ contains a set $A$
such that for every $\e>0$ the
$\e$-neighbourhoods of the translations of $A$ by finitely many
elements of $G$ always intersect. As a corollary,
amenability of $\pi$ is 
equivalent to the existence of a $G$-invariant mean on the
uniformly continuous bounded functions on $\s_\pi$.
As another corollary, a locally compact group $G$ is amenable if and
only if for every strongly continuous unitary representation of $G$ in an
infinite-dimensional Hilbert space $\mathcal H$ the system
$(\s_\H,G)$ has the property of concentration.
\end{abstract}

\maketitle

\setcounter{tocdepth}{1}
\tableofcontents

\section{Introduction} 
Let $\pi$ be a unitary representation 
of a group $G$ in a Hilbert space $\mathcal H$.
Then in particular $G$ acts on the unit sphere
$\s_{\mathcal H}$ of the space of representation.
The resulting topological dynamical system
$(\s_{\mathcal H},G,\pi)$ (which we will also denote by
$(\s_\pi,G)$) is thus a pretty common object in mathematics, 
and examining its properties from the dynamics viewpoint
could be worth while.

Dynamical systems of this kind have received plenty
of attention for finite-dimension\-al representations: 
in this case one can assume without loss in
generality that $G$ is a (compact)
Lie group, and actions of Lie groups on
finite-dimensional spheres are being studied intensely,
cf. e.g. \cite{GOn}. However,
for $\H$ infinite-dimensional much less 
appears to be known.
One obvious reason for that is the non-compactness of an 
infinite-dimensional sphere --- indeed, 
the main body of concepts and results in
present-day topological dynamics have substance for (locally) 
compact phase spaces only, cf. \cite{Aus,Ell,P2,dV1}. 
But what if one compactifies the dynamical system $(\s_{\mathcal H},G)$?

If we equip $\s_{\mathcal H}$ with the additive uniform structure
(determined by the norm), then $G$ acts on the sphere
by uniform isomorphisms. Denote by $\sigma\s_{\mathcal H}$ the 
maximal uniform, or Samuel,
compactification of the uniform space $\s_{\mathcal H}$,
that is, the Gelfand spectrum of the $C^\ast$-algebra of all bounded
uniformly continuous functions on the sphere.
Every uniform isomorphism of $\s_\H$ determines a unique
self-homeomorphism of
$\sigma\s_H$, and in this way $G$ acts on the compactum $\sigma\s_\H$.
Now we are in the realm of abstract topological dynamics,
where a traditional question of importance is: 
does a given compact $G$-flow
contain fixed points?

The existence of a fixed point in $\sigma\s_\H$ can be  
expressed in terms of the original system $(\s_\H,G)$. 
Following Milman \cite{M1,M2} and only
slightly extending a setting for his definition, 
let us introduce the following concept. 

\begin{defin}
Let $X=(X,{\mathcal U}_X)$ be a uniform space, and let
$F$ be a family of uniformly continuous self-maps of $X$. \par
A subset $A\subseteq X$ is called {\it essential} (for $F$) if
for every entourage of the diagonal 
$V\in{\mathcal U}_X$ and every finite collection of transformations
$f_1,f_2,\dots,f_n\in F$, $n\in\N$, one has
\begin{equation}
\bigcap_{i=1}^n f_iV[A]\neq\emptyset,
\label{milman}
\end{equation}
where 
$V[A]=\{x\in X\colon ~\mbox{for some $a\in A$,}~
(a,x)\in V\}$ is the $V$-neighbourhood of $A$. \par
One says that the pair $(X,F)$ has the {\it property of concentration}
if every finite cover $\gamma$ of $X$ contains an essential
set $A\in\gamma$. \tri
\label{milmandef}
\end{defin}

The property of concentration implies --- and, if $F$ is a group, is
equivalent to ---
the existence of a common fixed point for $F$ in $\sigma X$
(equivalently, in every $F$-equivariant uniform
compactification of the uniform space
$X$), cf. Proposition \ref{equiv} below. 

An important observation was made by Gromov and Milman 
\cite{GrM} (cf. also \cite{M1}): in a number of situations,
the concentration property of a dynamical system of the form
$(X,F)$, where $X$ is a uniform space infinite-dimensional
in some clear sense, is just another manifestation of
the phenomenon of concentration of measure on high-dimensional
structures \cite{FLM,Ma,M2,MS,St,Ta}. 
Among the results proved by Gromov and Milman 
\cite{GrM,M1} there are the following three. 

\begin{itemize}
\item If $G$ is 
abelian and $\dim\H=\infty$, then the pair  $(\s_{\mathcal H},G)$ has the 
property of concentration.

\item If $G$ is compact and $\H$ infinite-dimensional, 
then $(\s_{\mathcal H},G)$ has the 
property of concentration.

\item
The pair $(\s^\infty,U(\infty))$, where
$\s^\infty$ is the unit sphere of $l_2$  and 
$U(\infty)=\cup_{i=1}^\infty U(n)$,
has the property of concentration.
\end{itemize}

These results had led to the following natural
question, which, though not included by Gromov and Milman in their 
original paper \cite{GrM}, was later advertised by Milman in 
\cite{M1,M2}:
{\it does the pair $(\s_{\mathcal H},U({\mathcal H}))$
have the property of concentration for an infinite-dimensional Hilbert space
$\mathcal H$?}

The answer to the question is `No,' and a very simple
counter-example was constructed already in 1988 by Imre Leader 
\cite{L}, see Example \ref{leader} below. 
Unfortunately, the example was never published and, in particular,
the present author only learned about its existence after
his note \cite{P3} had appeared.
The existence of such an example suggests a deeper
reading of the above question: which groups
of unitary transformations have the concentration property on spheres
and which do not, and why? 
Some light on the point at issue was thrown by the present author in
\cite{P3}, where the following was proved.

\begin{itemize}
\item A discrete group $G$ is amenable if and only if the dynamical system
$(\s_\H,G,\pi)$ has the property of concentration for every unitary
representation $\pi$ of $G$ in an infinite-dimensional Hilbert space
$\H$. (Equivalently: for the left regular
representation $\pi_2$.)
\end{itemize}

It is evident that if $H$ is a subgroup of $G$ and $(\s_\H,G)$ has the
property of concentration, then so does $(\s_\H,H)$. Therefore,
for example, the pair formed by the unit sphere of the space
$L_2(F_2)$ and the full unitary group of this space does not have the
property of concentration, where $F_2$ denotes the free group on two
generators.

It remained still unclear whether or not the above result could be extended to
locally compact groups. Indeed attempts to link properties of the
$G$-flow $(\sigma\s_\H,G,\pi)$ with topologo-algebraic properties 
of a non-discrete group $G$ 
encounter the following difficulty: in general, the
extended action of $G$ on the Samuel compactification
$\sigma\s_\H$ is no longer continuous, and thus the dynamical system
$(\sigma\s_\H,G,\pi)$ does not even `remember' the topology of $G$.

This observation suggests that the concentration property of the
system $(\s_\H,G,\pi)$ has to do not with the amenability of the acting
group $G$ as such, but rather with the amenability of the representation
$\pi$ as defined by Bekka \cite{B}. Adopting this viewpoint
(suggested by Pierre de la Harpe after he got acquainted with our e-print
\cite{P3}) turns out to be very fruitful. 

\begin{defin}
According to Bekka \cite{B}, 
a unitary representation $\pi$ of a group $G$ in a Hilbert space
$\H$ is {\it amenable} if there exists a $G$-invariant
state $\phi$ on the algebra of bounded operators ${\mathcal L}(\H)$.
It means that $\phi\in{\mathcal L}(\H)^\ast$, $\phi\geq 0$,
$\phi({\mathbb I})=1$, and $\phi(\pi(g)T\pi(g)^{-1})=\phi(T)$ for
every $T\in{\mathcal L}(\H)$ and every $g\in G$. Such a functional
$\phi$ is called a $G${\it -invariant mean}. 
\tri\end{defin}

This concept unifies several previous
theories of amenability, and in particular a locally compact group $G$
is amenable if and only if every strongly continuous representation of
$G$ is amenable \cite{B}. Notice also that amenability of a unitary
representation does not depend on the topology of $G$.

Amenability of a representation turns out to be a necessary 
prerequisite for the concentration property.
\vskip .3cm

\noindent{\bf Corollary \ref{if}.} {\it
Let $\pi$ be a unitary representation of a group $G$ in a
Hilbert space $\H$. If the dynamical system 
$(\s_\H,G,\pi)$ has the concentration property,
then the representation $\pi$ is amenable.} \qed
\vskip .3cm

We are of course more interested in trying to reverse this statement
{\it \`a la} Gromov and Milman. As Example \ref{ex2} shows, 
the $G$-flow 
$(\s_\H,G,\pi)$ need not have the concentration property even
if a representation $\pi$ of a group $G$ in an infinite-dimensional
space $\H$ is amenable. Nevertheless, excluding `trivially
amenable' representations leads to the following result.
\vskip .3cm

\noindent{\bf Theorem \ref{conc}.} {\it
 Let $\pi$ be a unitary representation of a group $G$ in
an infinite-dimensional Hilbert space $\H$. If every
subrepresentation of $\pi$ having finite codimension is amenable
(that is, $\pi$ is not of the form
$\pi_1\oplus\pi_2$, where $\pi_1$ is finite-dimensional and
$\pi_2$ is non-amenable), then 
the dynamical system $(\s_\H,G,\pi)$ has the concentration property.}
\qed
\vskip .3cm

The proof is again based on the technique of concentration of measure on
high-dimensional structures.

Now we are able to derive a number of definitive results linking
the amenability property of a unitary representation $\pi$ with the
concentration property of the system $(\s_\H,G,\pi)$.
We begin with a description of subgroups of the full unitary
group $U(\H)$ whose action on the unit sphere has the concentration
property.
\vskip .3cm

\noindent
{\bf Theorem \ref{either}.} {\it
 Let $\pi$ be a unitary representation of a group
$G$ in a Hilbert space $\H$. The system
$(\s_\H,G,\pi)$ has the concentration property if and only if
\begin{itemize}
\item either $\pi$ has a non-zero invariant vector, or
\item $\dim\H=\infty$ and
every subrepresentation of $\pi$ having finite codimension
is amen\-able. \qed
\end{itemize}
}
\vskip .3cm

We can now extend our criterion of amenability
from discrete groups \cite{P3} to all locally compact ones.
\vskip .3cm

\noindent
{\bf Theorem \ref{charact}.}
{\it 
A locally compact group $G$ is amenable if and only if for every
strongly continuous unitary representation $\pi$ of $G$ in an
infinite-dimensional Hilbert space $\H$, the dynamical system
$(\s_\H,G,\pi)$ has the concentration property.} \qed
\vskip .3cm

In their turn, amenable representations can be characterized
in terms of the concentration property.
\vskip .3cm

\noindent
{\bf Theorem \ref{amenability}.} {\it 
A unitary representation $\pi$ of a group $G$ in a Hilbert space $\H$
is amenable if and only if 
\begin{itemize}
\item 
either $\pi$ contains a finite-dimensional
subrepresentation, or 
\item the $G$-space $(\s_\H,G,\pi)$ has the
concentration property. \qed
\end{itemize}}
\vskip .3cm

One of the applications is to the `L\'evy-type integral'
for functions on the unit sphere in an infinite-dimensional Hilbert space.
(Cf. \cite{Gri, M2}.)
\vskip .3cm

\noindent
{\bf Theorem \ref{levy}.} {\it 
Let $\pi$ be a unitary representation of a group
$G$ in a Hilbert space $\H$. The following conditions are
equivalent.
\begin{itemize}
\item The space $C^b_u(\s_\H)$ of all 
bounded uniformly
continuous functions on the unit sphere $\s_\H$  
admits a $G$-invariant mean.
\item The representation $\pi$ is amenable. \qed
\end{itemize}}
\vskip .3cm

Also in Section \ref{dynamical} 
we establish a number of dynamical corollaries
listed in our C.r. note \cite{P3}
without proofs.

\section{\label{fixed}Fixed points and concentration property} 

Let $X=(X,{\mathcal U}_X)$ be a uniform space.
The {\it Samuel compactification}, or else
the {\it maximal uniform compactification,} of $X$ \cite{Eng}
is a Hausdorff compact space $\sigma X$
together with a uniformly continuous mapping $i_X\colon X\to\sigma X$
such that every uniformly continuous mapping $f$ of $X$ to an
arbitrary compact Hausdorff space $K$ factors through $i_X$, that is,
there exists a continuous mapping $\bar f\colon \sigma X\to K$ 
with $f=\bar f \circ i_X$. (Recall that every compact space supports
a unique compatible uniform structure.) In particular, it follows easily
that the image $i_X(X)$ is everywhere dense in $\sigma X$.

The Samuel compactification $\sigma X$ is
the completion of the uniform space $(X,{\mathcal C}^\ast(X))$,
where ${\mathcal C}^\ast(X)$ is the finest totally bounded uniform
structure on $X$ contained in ${\mathcal U}_X$. The uniformity
${\mathcal C}^\ast(X)$ is at the same time 
the coarsest uniformity making each
bounded uniformly continuous complex-valued 
function on $(X,{\mathcal U}_X)$ 
uniformly continuous on $(X,{\mathcal C}^\ast(X))$. 

The Stone-\u Cech
compactification, $\beta X$, of a Tychonoff topological space $X$ is
a special case of the Samuel compactification recovered
if $X$ is equipped with the finest compatible uniformity.

In particular, every uniformly continuous mapping
$f\colon X\to X$ determines
a unique continuous mapping $\bar f\colon \sigma X\to \sigma X$.
If $f$ is a uniform automorphism of $X$, then 
$\bar f$ is a self-homeomorphism of $\sigma X$.

Let $F$ be a family of uniformly continuous self-maps of a uniform space
$X$. A compactification $(K,i)$ of $X$ (that is, a pair formed by a
compact space $K$ and a uniformly continuous mapping $i\colon X\to K$
with an everywhere dense image) is called $F${\it -equivariant} if
for each $f\in F$ there exists a (necessarily unique)
continuous mapping $\tilde f\colon K\to K$
satisfying $\tilde f\circ i= i\circ f$. It follows that the
Samuel compactification of $X$ is $F$-equivariant.

The Samuel compactification can also be described as the
Gelfand spectrum of the commutative $C^\ast$-algebra
$C^b_u(X)\cong C(\sigma X)$
of all bounded uniformly
continuous complex-valued functions on a uniform space $X$ equipped
with the supremum norm. Since every uniformly continuous 
mapping $f\colon X\to X$
gives rise to a unital $C^\ast$-algebra endomorphism 
$f^\ast$ of $C^b_u(X)$, one can talk of $F$-invariant
means (in the $C^\ast$-algebraic terminology, states) on $C^b_u(X)$.

\begin{prop}
\label{equiv} For a family $F$ of automorphisms
of a uniform space $X=(X,{\mathcal U}_X)$
the following are equivalent.
\begin{enumerate}
\item\label{one} The pair $(X,F)$ has the property of concentration.
\item\label{six} For every finite subfamily $F_1\subseteq F$, the pair
$(X,F_1)$ has the property of concentration.
\item \label{two} The pair $(\sigma X,F)$ has the property of concentration.
\item\label{three} The family $F$ has a common fixed point in the
Samuel compactification of $X$.
\item\label{five} The family $F$ has a common fixed point in every
$F$-equivariant uniform compactification of $X$.
\item\label{four} There exists an $F$-invariant multiplicative mean on  
the space $C_u^b(X)$.
\end{enumerate}
If $F$ is a family of uniformly continuous self-maps of $X$, then
{\rm (\ref{one}) $\Leftrightarrow$ (\ref{six}) 
$\Rightarrow$ (\ref{two}) $\Leftrightarrow$
(\ref{three}) $\Leftrightarrow$ (\ref{five})
$\Leftrightarrow$ (\ref{four}).}
\end{prop}

\begin{proof} 
(\ref{one}) $\Rightarrow$ (\ref{six}): obvious.

(\ref{six}) $\Rightarrow$ (\ref{one}): 
Let $\gamma$ be a finite cover of $X$.
For every finite $F_1\subseteq F$ denote
by $\gamma_{F_1}$ the (non-empty, finite) collection of
all $F_1$-essential elements of $\gamma$. Clearly, whenever $F_1\subseteq F_2$,
one must have $\gamma_{F_2}\subseteq\gamma_{F_1}$. The compactness 
(or rather finiteness)
considerations lead one to conclude that 
\begin{equation}
\bigcap_{F_1\subseteq F,~
\vert F_1\vert<\infty}\gamma_{F_1}\neq\emptyset,
\end{equation}
 thus finishing the proof: every element $A$ of the above intersection
is $F$-essential.

(\ref{one}) $\Rightarrow$ (\ref{two}): if $\gamma$ is a
finite cover of $\sigma X$, then at least one of the sets $A\cap X$,
$A\in\gamma$, is $F$-essential in $X$, and it follows 
that $A$ is $F$-essential in $\sigma X$.

(\ref{two}) $\Rightarrow$ (\ref{three}): 
 emulates a proof of
Proposition 4.1 and Theorem 4.2 in \cite{M1}. 

There exists a point $x^\ast\in \sigma X$ whose every neighbourhood is
essential: assuming the contrary, one can cover the compact space
$\sigma X$ with open $F$-inessential sets 
and select a finite subcover containing
no $F$-essential sets, a contradiction.

We claim that $x^\ast$ is a common fixed point for $F$.
Assume it is not so.
Then for some $f\in F$ one has $fx^\ast\neq x^\ast$.
Choose an entourage, $W$, of the unique compatible 
uniform structure on $\sigma X$ with the property
$W^2[x^\ast]\cap W^2[fx^\ast]=\emptyset$.
Since $f$ is uniformly continuous, there is a $W_1\subseteq W$ with
$(x,y)\in W_1\Rightarrow (fx,fy)\in W$ for all $x,y$.
Since $f(W_1[x^\ast])\subseteq W[f(x^\ast)]$,
 we conclude that
$W_1[x^\ast]$ is $F$-inessential (with $V=W$), a contradiction. 

(\ref{three}) $\Rightarrow$ (\ref{five}): Let 
$(K,i)$ be an $F$-equivariant
compactification of $X$. There exists a unique continuous 
$j\colon \sigma X\to K$ with $j\circ i_X=i$, and it follows easily
that $j$ is an $F$-equivariant mapping, that is, 
$j\circ \bar f= \tilde f\circ j$ for every $f\in F$.
The image of an $F$-fixed point $x^\ast$ under $j$ is an $F$-fixed
point in $K$.

(\ref{five}) $\Rightarrow$ (\ref{three}): trivial.

(\ref{three}) $\Leftrightarrow$ (\ref{four}): 
fixed points in the Gelfand space $\sigma X$ of
the commutative $C^\ast$-algebra $C_u^b(X)$ 
correspond to $F$-invariant multiplicative means
(states) on $C_u^b(X)$.

(\ref{three}) $\Rightarrow$ (\ref{two}):
If $\gamma$ is a finite cover of $\sigma X$, then there is an $A\in\gamma$ 
containing an $F$-fixed point, and such an $A$ is clearly $F$-essential. 

(\ref{three}) $\Rightarrow$ (\ref{six}): this is the only implication
where we assume $F$ to be a group.

Without loss in generality and replacing $X$ with its separated reflection
if necessary, one can assume that $X$ is a separated uniform space
(that is, $\cap{\mathcal U}=\Delta_X$): indeed, the Samuel compactifications
of a uniform space $X$ and of its separated
reflection are canonically homeomorphic.
Thus, we can identify $X$ (as a topological, not uniform space!)  
with an everywhere dense subspace of $\sigma X$.

If now $\gamma$ is a finite cover of $X$, then the closures of all 
$A\in\gamma$ taken in $\sigma X$ cover the latter space, and so
there is an $A\in\gamma$ with 
$\operatorname{cl}\,_{\sigma X}(A)$ containing
an $F$-fixed point $x^\ast\in\sigma X$. 
We claim that $A$ is $F$-essential in $X$. 

To prove this, we need a simple fact of general topology. Let
$B_1,\dots,B_n$ be subsets of a uniform space $X$ satisfying the
condition $V[B_1]\cap \dots\cap V[B_n]=\emptyset$ for some entourage
$V\in {\mathcal U}_X$. Then the closures of $B_i$, $i=1,2,\dots,n$
in the Samuel compactification $\sigma X$ have no point in common:
$\operatorname{cl}_{\sigma X}(B_1)\cap \dots \cap
\operatorname{cl}_{\sigma X}(B_n)=\emptyset$.

[Let $\rho$ be a uniformly
continuous bounded pseudometric on $X$ subordinated to the entourage
$V$ in the sense that $(x,y)\in V$ whenever $x,y\in X$ and
$\rho(x,y)<1$.
For each $i=1,\dots,n$ and $x\in X$ 
set $d_i(x)=\inf\{\rho(b,x)\colon b\in B_i\}$. 
The real-valued functions $d_i$
are uniformly continuous (indeed $\rho$-Lipschitz-1) 
and bounded on $X$, and therefore extend to
(unique) continuous functions $\tilde d_i$ on $\sigma X$. 
If there existed a common point, $x^\ast$, for the
closures of all $B_i$ in $\sigma X$, then all $\tilde d_i$ would 
vanish at $x^\ast$ and consequently  for any given $\e>0$ there 
would exist an $x\in X$ with $d_i(x)<\e$, and in particular
$x\in V[B_i]$, for all $i$.
However, for every $x\in X$, there is an $i$ with $x\notin V[B_i]$.]

Now assume that $\cap_{f\in F_1}f(V[A])=\emptyset$ for 
some $V\in{\mathcal U}_X$,
where $F_1$ is a finite subfamily of $F$. Every $f\in F_1$ 
has a uniformly continuous inverse
$f^{-1}\in F$, and there is an entourage $V_1\in {\mathcal U}_X$
with $(x,y)\in V_1\Rightarrow (f^{-1}x,f^{-1}y)\in V$ 
for all $x,y\in X$ and $f\in F_1$. Let $x\in V_1[f(A)]$, that is,
$(x,f(a))\in V_1$ for some $a\in A$. Then 
$(f^{-1}(x),a)\in V$, that is, $y=f^{-1}(x)\in V[A]$,
and consequently $x=f(y)\in f(V[A])$. We conclude:
$V_1[f(A)]\subseteq f(V[A])$ for each $f\in F_1$, and therefore
$\cap_{f\in F_1}V_1[f(A)]=\emptyset$ as well.

The above observation from uniform topology implies that 
$\cap_{f\in F}\operatorname{cl}_{\sigma X}(f(A))=\emptyset$.
Since extensions of $f$ to $\sigma X$ are homeomorphisms,
$\operatorname{cl}_{\sigma X}(f(A))=
\bar f(\operatorname{cl}_{\sigma X}(A))$, and consequently
$\cap_{f\in F_1}f(\operatorname{cl}_{\sigma X}(A))=\emptyset$,
a contradiction because the intersection contains $x^\ast$.
\end{proof}

\begin{remark} As pointed out in \cite{Gl} on a similar
occasion, the condition of uniform equicontinuity of $F$, imposed in
\cite{GrM,M1,M2}, is superfluous. \tri
\end{remark}

\begin{example} Here is a simple example showing that in general
$(\ref{three})\not\Rightarrow (\ref{six})$ if $F$ consists just of
uniformly continuous mappings. Let $X=(-1,0)\cup (0,1)$ with the
additive uniform structure, and set $F=\{f_1,f_2\}$, where
$f_1(x)=\abs x$ and $f_2(x)=-\abs x$. Then $\sigma X=[-1,1]$ and
$x^\ast=0$ is an $F$-fixed point in $X$. At the same time, the pair
$(X,F)$ does not have the property of concentration, as can be seen
from considering the cover $\gamma=\{(-1,0),(0,1)\}$, both elements of
which are $F$-inessential.
\tri\end{example}

\begin{remark} If one however replaces the condition (\ref{milman}) in
the Definition \ref{milmandef} of the concentration property
with the condition
\begin{equation}
\bigcap_{i=1}^n V[f_i(A)]\neq\emptyset,
\end{equation}
then all six conditions in Proposition \ref{equiv} become equivalent for
an arbitrary family $F$ of uniformly continuous self-maps of $X$.
Perhaps, using the concept of concentration property so modified is a
sensible thing to do. \tri
\end{remark}

\begin{example}
\label{almost}
A unitary representation
$\pi$ of a group $G$ {\it has almost invariant vectors}
if for every finite $F\subseteq G$ and every $\e>0$ there is a
$\xi$ in the space of representation $\H_\xi$ with $\norm\xi=1$
and $\norm{g\cdot\xi-\xi}<\e$ for every $g\in F$. 
As can be easily seen at the level of definitions, 
if a representation $\pi$ has almost invariant vectors, then the
system $(\s_\H,G,\pi)$ has the concentration property.

The converse is not true. The simplest example possible is the
representation of $\Z_2\cong\{1,-1\}$ in $l_2$ by scalar multiplication.
Non-existence of almost invariant vectors is manifest,
yet according to the results by Gromov and Milman
cited in the Introduction, the system $(\s_\infty,\Z_2)$ has the
concentration property.
\tri
\end{example}

The following result will be employed in Sections \ref{fromam}
and \ref{vis}.

\begin{prop}
Let $\pi_i$, $i=1,2$, be unitary representations of a group $G$ in
Hilbert spaces $\H_i$.
Let $\pi=\pi_1\oplus\pi_2$ be the direct sum representation. 
Then the following are equivalent.
\begin{enumerate}
\item The dynamical system $(\s_{\H_1\oplus\H_2},G,\pi_1\oplus\pi_2)$ 
has the concentration property.
\item  At least one of the systems $(\s_{\H_i},G,\pi_i)$,
$i=1,2$ has the concentration property.
\end{enumerate}
\label{sum}
\end{prop}

\begin{proof} Since the unit spheres $\s_i=\s_{\H_i}$, $i=1,2$, are
contained in $\s_\H$ in a canonical way both as uniform
subspaces and $G$-subspaces, 
it follows that the compactifications $\sigma\s_i$, $i=1,2$, are compact
$G$-subflows of $\sigma\s_\H$. This establishes
(2) $\Rightarrow$ (1).

Now assume that both systems  $(\s_{\H_i},G,\pi_i)$,
$i=1,2$ do not have the concentration property. Then there exist
finite covers
$\gamma_j$  of $\s_j$, $j=1,2$, a finite collection $g_1,\dots,g_n$ of
elements of $G$, and an $\e>0$ having the property that for
every $A\in\gamma_j$, $j=1,2$:
\[\cap_{i=1}^n\pi_{j}(g_i)({\mathcal O}_\e(A))=\emptyset.\]
(Here and below it is convenient to assume that the $\e$-neighbourhood of
$A$ in the sphere, 
${\mathcal O}_\e(A)$, is formed with respect to the geodesic distance.)

For every $j=1,2$ and each $A\in\gamma_j$ set
\[\tilde A=\left\{x\in\s_\H\colon \norm{\pi_jx}\geq\frac{\sqrt 2}{2}
~~\mbox{and}~~ \frac{\pi_jx}{\norm{\pi_jx}}\in A \right\}.\]
The collection $\cup_{j=1,2}\{\tilde A\colon A\in\gamma_j\}$ 
covers $\s_\H$, and it is 
easy to see that for each $j=1,2$ and each $A\in\gamma_j$ one has
\[\cap_{i=1}^n\pi_{j}(g_i)({\mathcal O}_\delta(\tilde A))=\emptyset\]
whenever $\delta$ is sufficiently small, for example
$\delta< \min\{\e/3, \pi/8\}$. 
\end{proof}

\section{\label{fromconc}From concentration property to amenability}

\begin{prop} Let $\pi$ be a representation of a group $G$ in a Hilbert
space $\H$. If the space $C^b_u(\s_\H)$ of all bounded uniformly
continuous functions on the unit sphere $\s_\H$ admits a $G$-invariant
mean, then the representation $\pi$ is amenable.
\label{inv}
\end{prop}

\begin{proof} Let $\psi\colon C^b_u(\s_\H)\to\C$ denote a $G$-invariant
mean on the unit sphere $\s_\H$, that is, 
a positive functional of norm $1$, taking the
function $1$ to $1$ and such that $\psi(_gf)=\psi(f)$ for all
$g\in G$ and all $f\in C^b_u(\s_\H)$, where $_gf(x):=f(\pi(g)x)$.
For every bounded linear operator $T$ on
$\H$ define a function $f_T\colon\s_\H\to\C$ by
\begin{equation}
\s_\H\ni\xi\mapsto f_T(\xi):=(T\xi, \xi)\in \C.
\end{equation}
Then $f_T$ is bounded (by $\norm T$) and Lipschitz with constant
$2\norm T$:
\begin{eqnarray}
f_T(\xi)-f_T(\eta)&=&
(T\xi,\xi)-(T\eta,\eta)\nonumber \\
&=&
(T\xi,\xi-\eta) + (T(\xi-\eta),\eta)\nonumber \\
&\leq& 2\norm T\cdot\norm{\xi-\eta}.
\end{eqnarray}
Therefore, $f_T\in C^b_u(\s_\H)$. Set 
\begin{equation}
\phi(T):=\psi(f_T).
\end{equation}
The following properties of the mapping $\phi\colon L(\H)\to\C$ are
obvious.
\smallskip

\begin{enumerate}
\item $\phi$ is linear. [If $T,S\in L(\H)$ and $\lambda,\mu\in\C$, then
for every $\xi\in\s_\H$ one has $f_{\lambda T+\mu S}(\xi)=
((\lambda T+\mu S)\xi,\xi)=\lambda 
(T\xi,\xi)+\mu(S\xi, \xi)=\lambda f_T(\xi)+\mu f_S(\xi)$, 
and consequently
$\phi(\lambda T+\mu S)=\psi(f_{\lambda T+\mu S})=
\lambda \psi( f_T)+\mu
\psi( f_S)=\lambda \phi(T)+\mu\phi(S)$.]
\smallskip

\item $\phi$ is positive. [If $T\geq 0$, then $f_T(\xi)=(T\xi,\xi)\geq 0$
for all $\xi$, therefore $\phi(T)=\psi( f_T)\geq 0$.]
\smallskip

\item $\phi({\mathbb I})=1$. [$f_{\mathbb I}(\xi)=(\xi,\xi)\equiv 1$.]
\smallskip

\item $\phi$ is $G$-invariant.
[Let $g\in G$ and $T\in L(\H)$. Then 
\begin{eqnarray}
f_{\pi(g)T\pi(g)^{-1}}(\xi)&=&
(\pi(g)T\pi(g)^{-1}(\xi),\xi)\nonumber \\
&=&(T\pi(g)^{-1}\xi,\pi(g)^{-1}\xi)\nonumber \\
&=& f_T(\pi(g)^{-1}\xi )\nonumber \\
&=&
_{g^{-1}}(f_T)(\xi).
\end{eqnarray}
Therefore,
\begin{eqnarray}
\phi(\pi(g)T\pi(g)^{-1}) &=& \psi( f_{\pi(g)T\pi(g)^{-1}})
\nonumber \\
&=& \psi(_{g^{-1}} f_T) \nonumber \\
&=& \psi(f_T) \nonumber \\
&=& \phi(T).]
\end{eqnarray}
\end{enumerate}

The conditions (1)-(3) imply that 
$\phi$ is also bounded of norm $1$ (cf. \cite{Sa}, Prop.
1.5.1), and by (4) $\phi$ is a $G$-invariant
mean on $L(\H)$, as required.
\end{proof}

\begin{remark} The above result will be inverted in the concluding Section
(Thm. \ref{levy}), leading to a new equivalent definition of amenable
representations, very much in the classical spirit. \tri
\end{remark}

\begin{corol} 
\label{if}
Let $\pi$ be a representation of a group $G$ in a Hilbert
space $\H$. If the dynamical system $(\s_\H,\pi,G)$ has the
concentration property, then the representation $\pi$ is amenable.
\end{corol}

\begin{proof} 
Let $x^\ast$ be a $G$-fixed point in the Samuel compactification
$\sigma\s_\H$. For every $f\in C^b_u(\s_\H)$ set $\psi(f)=\tilde f(x^\ast)$,
where $\tilde f$ is the unique continuous extention of $f$
to $\sigma\s_\H$. Clearly, $\psi$ is a $G$-invariant mean on
$C^b_u(\s_\H)$.
\end{proof}

\begin{remark}
The converse statement is false: indeed, every finite-dimensional
representation is amenable \cite{B}, but for such a representation 
the concentration property is equivalent to the existence of a non-zero
invariant vector.

As can be seen from results of Gromov and Milman, the assumption of
infinite-dimensionality of the space of representation is essential
for deriving the concentration property of the sphere. This
observation will be reinforced in the subsequent sections of our 
article. Strictly speaking, Corollary \ref{if} cannot be inverted 
even for infinite-dimensional
Hilbert spaces (Example \ref{ex2} below). However, if one dismisses
`trivially amenable' representations (that is, those whose amenability
stems from the existence of a single finite-dimensional
representation), then the concentration property of the sphere is back, 
cf. Theorem \ref{conc}. \tri
\end{remark}

\begin{corol}
Let $G$ be a locally compact group. Denote by $\pi_2$ the left
regular representation of $G$, and by $\s_2$ the unit sphere in
the space $L_2(G)$. If the system
$(\s_2,G,\pi_2)$ has the concentration property, then $G$ is
an amenable LC group.
\label{l2}
\end{corol}

\begin{proof} It is enough to apply the following result of Bekka
(\cite{B}, Thm. 2.2): the left regular representation $\pi_2$ of a LC group
$G$ is amenable if and only if $G$ is amenable.

It is instructive to look at the direct proof as well.
Let $x^\ast$ be a $G$-fixed point in $\sigma\s_2$.
For every Borel set $A\subseteq G$ and each 
$f\in\s_2$ set $z_A(f)=\norm{\chi_A\cdot f}^2$, where $\chi_A$ denote,
as usual, the characteristic function of $A$, and the dot stands for
the muptiplication of (equivalence classes of) functions.
Since the mapping $f\mapsto\norm f^2$ is $2$-Lipschitz on $\s_2$,
so is the function $z_A\colon \s_2\to\R$. 
Being also bounded, $z_A\in C^b_u(\s_2)$. Denote by $\tilde{z_A}$
the unique continuous extension of $z_A$ to the
Samuel compactification of the sphere, and set 
$m(A)=\tilde{z_A}(x^\ast)$. Then $m$ is a finitely-additive left-invariant
normalized measure on Borel subsets of $G$, vanishing on locally null
sets, and consequently $G$ is amenable.
\end{proof}

\begin{corol}  The unit sphere of an infinite-dimensional
Hilbert space $\mathcal H$ admits no left invariant means 
on bounded uniformly continuous functions with respect to the full unitary 
group $U({\mathcal H})$.
\label{means}
\end{corol}

\begin{proof}
It is enough to make an obvious remark:
an $U({\mathcal H})$-invariant mean on $C^b_u(\s_{\mathcal H})$
is invariant with respect to the action of every group
$G$ represented in $\mathcal H$ by unitary operators. 
Since $\mathcal H$ is infinite-dimensional,
one can find a non-amenable discrete group $G$ of the same cardinality
as is the density character of $\mathcal H$, and to realize $\mathcal H$ as $l_2(G)$.
(For example, take as $G$
the free group of rank equal to the density character of
$\mathcal H$.)
Now one can apply Prop. \ref{inv}.
\end{proof}

\begin{remark} 
The above result means, in essence, that 
there exists no L\'evy-type integral of uniformly continuous functions
on the sphere that is invariant under the
action of the full unitary group. (Cf. \cite{Gri, M2}.) Soon we will see
(Theorem \ref{levy})
that the existence of a $G$-invariant L\'evy-type integral on the unit
sphere $\s_\H$ is in fact
equivalent to the amenability of the representation of $G$ in $\H$. \tri
\end{remark}

\begin{corol} 
\label{hilb}
Let $\mathcal H$ be a Hilbert space. The pair
$(\s_{\mathcal H},U({\mathcal H}))$ does not have the concentration property.
\end{corol}

\begin{proof}
 If $\mathcal H$ 
is infinite-dimensional, the statement follows from 
Corol. \ref{l2} or \ref{means}. 
If $\dim{\mathcal H}<\infty$, the unitary group $U(n)$ possesses no
non-zero invariant vectors, and there is no concentration property
in a trivial way.
\end{proof}

\begin{example}
\label{leader} 
The first counter-example showing the absence of 
the concentration property
of the system $(\s^\infty,U(l_2))$ was obtained in 1988 by
Imre Leader \cite{L} (unpublished). This example is remarkably simple ---
indeed, as simple as one can probably ever get.
Here we reproduce it upon a kind permission from
the author. 

Let $\s^\infty$ denote the unit sphere in the space
$l_2(\N)$. For a $\Gamma\subset\N$, let $p_X$ stand for the orthogonal
projection of $l_2(\N)$ onto its subspace $l_2(\Gamma)$.
 Denote by $E$ (respectively, $F$) the set of all even
(resp., odd) natural numbers, and let
\[A=\left\{x\in \s^\infty \colon \norm{P_Ex}
\geq {\sqrt 2}/2\right\},\]
\[B=\left\{x\in \s^\infty \colon \norm{P_Ex}
\leq {\sqrt 2}/2\right\}.\]
Clearly, $A\cup B=\s^\infty$. At the same time, both $A$ and $B$ are
inessential. To see this, let $E_1,E_2,E_3$ be three
arbitrary disjoint infinite
subsets of $\N$, and let $\phi_i\colon \N\to\N$ be bijections
with $\phi_i(E)=E_i$, $i=1,2,3$. Let $g_i$ denote the unitary
operator on $l_2(\N)$ induced by $\phi_i$. 
Now
\[g_i(A)=\left\{x\in\s^\infty\colon \norm{p_{E_i}x}\geq
{\sqrt 2} /2\right\},\]
and consequently
\[{\mathcal O}_\e(g_i(A))\subseteq 
\left\{x\in\s^\infty\colon \norm{p_{E_i}x}\geq
({\sqrt 2} /2)-\e\right\}.\]
Thus, as long as $\e<{\sqrt 2}/2-{\sqrt 3}/3$,
we have
\[\cap_{i=1}^3 {\mathcal O}_\e(g_i(A))=\emptyset.\]
The set $B$ is treated in exactly the same fashion.
\qed
\end{example}

Another such counter-example can be obtained through
combining the proof of Corol. \ref{l2} with von Neumann's proof of
non-amenability of $F_2$ (cf. e.g. \cite{Pa}, ex. 0.6).
In this form the link between amenability and concentration property
becomes very transparent, while the construction remains remarkably
similar to that in Leader's example \ref{leader}.

\begin{example} \cite{P3}
\label{expl} 
Let $a,b$ be free generators of $F_2$, and let 
$\pi=\pi_2$ be the left regular representation of $F_2$ in
${\mathcal H}=l_2(F_2)$; we will write $xf$ for $\pi_x(f)$.
Denote by $W_n$ the collection of all words whose irreducible
representation starts with $a^n$, $n\in\Z$.
Set 
\[A_1=\{f\in\s_2\colon \norm{\chi_{W_0}\cdot f}\leq 1/3\},\]
\[A_2=\{f\in\s_2\colon \norm{\chi_{W_0}\cdot f}\geq 1/3\},\]
 and
\[F=\{a,a^2,a^3,a^4,b\}.\]
Clearly, 
$\s_2=A_1\cup A_2$. Both $A_1$ and $A_2$ are $F$-inessential.
Indeed, if $f\in A_1$, then $\norm{\chi_{W_0}\cdot bf}\geq
\norm{\chi_{F_2\setminus W_0}\cdot f}\geq 2/3$
and consequently 
\[{\mathcal O}_{1/12}(A_1)\cap {\mathcal O}_{1/12}(bA_1)
=\emptyset.\]
If $f\in A_2$, there is an $i\in \{1,2,3,4\}$ such that
$\norm{\chi_{W_{-i}}\cdot f}<1/6$, and consequently 
$\norm{\chi_{W_0}\cdot a^if}<1/6$, meaning that
\[{\mathcal O}_{1/12}(A_2)\cap \cap_{i=1}^4{\mathcal O}_{1/12}(a^iA_2)
=\emptyset.\]
\qed
\end{example}

\section{\label{dynamical}Dynamical corollaries}
Let $G$ be a topological group and let $X$ be a topological
$G$-space, that is, a topological space equipped with
a {\it continuous} action of $G$.
The {\it maximal $G$-compactification} of 
$X$ is a compact $G$-space $\alpha_G(X)$ together with a
morphism of $G$-spaces (that is, a $G$-equivariant continuous
mapping) $i\colon X\to \alpha_G(X)$ such that any morphism
from $X$ to a compact $G$-space uniquely 
factors through $i$ \cite{dV,dV1,Meg}. Necessarily, the image $i(X)$ is
everywhere dense in $\alpha_G(X)$, though, somewhat surprisingly, 
the mapping
$i$ need not be a homeomorphic embedding --- in fact, it can be
even a constant mapping
for a nontrivial $G$-space $X$, cf. \cite{Meg}.

By $G/H_\Rsh$ we denote the left factor-space $G/H$ of a topological
group $G$ by a closed subgroup $H$, equipped with the uniformity
whose basis is formed by entourages of the form
\[V_\Rsh=\{(xH,yH)\colon xy^{-1}\in V\},\]
where $V$ is a neighbourhood
of $e_G$. One can show, using results from \cite{Meg}, that in
general the uniform space $G/H_\Rsh$ need not be separated and
can even induce the indiscrete topology.

If $H=\{e_G\}$, then we obtain the uniform space $G_\Rsh$, which is
always separated and induces the topology of the group $G$.

We say that a function $f\colon G/H\to\R$ is 
$\Rsh$-uniformly continuous
($\Rsh$-u.c.), if it satisfies the condition:
for every $\e>0$, there is a $V\ni e$ with
\begin{equation}
\forall x,y\in G,~~ xy^{-1}\in V ~~\Rightarrow ~~\vert f(xH)-f(yH)\vert<\e
\end{equation}

\begin{remark}
We deliberately avoid using the `right/left uniformly continuous'
terminology, because the mathematical
community seems to be divided into two groups of roughly the same size,
one of them calling the $\Rsh$-u.c. functions `right'
uniformly continuous, the other `left' uniformly continuous; references
to the both kinds of usage are given in \cite{P1}.

Our system of notation, suggested in \cite{P1,P2}, has a 
mnemonic advantage:
the symbol $\Rsh$ (in \TeX, $\$\backslash${\tt Rsh}$\$$)
reminds of the position of the
inversion symbol in the expression $xy^{-1}$. The functions satisfying
the property 
\begin{equation}
\forall x,y\in G,~~ x^{-1}y\in V ~~\Rightarrow ~~\vert f(xH)-f(yH)\vert<\e
\end{equation}
are naturally called $\Lsh$-uniformly continuous.
\tri
\end{remark}

Notice that bounded
$\Rsh$-uniformly continuous functions on $G/H$ 
are identified in an obvious way with bounded
$\Rsh$-u.c. functions on $G$ that are constant
on each left coset $xH$, $x\in G$. Their totality forms a $G$-invariant
$C^\ast$-subalgebra of $C^b_u(G_\Rsh)$, which we will denote 
by $C^b_u(G/H_\Rsh)$. 

The following result must be known, but it is difficult to find
an exact reference.

\begin{prop}
The maximal $G$-compactification of the left 
topological $G$-space $G/H$ coincides with
the Samuel compactification of $G/H_\Rsh$. 
\label{max}
\end{prop}

\begin{proof}
Since the left regular representation
of $G$ in $C^b_u(G_\Rsh)$ (defined by $(gf)(x)=f(gx)$)
 is well known (and easily checked) to be
strongly continuous \cite{Te,dV,Aus,P2}, so is the
subrepresentation of $G$ in $C^b_u(G/H_\Rsh)$. Now it follows from
a result of Teleman \cite{Te} that the action of $G$ on the Gelfand
spectrum of $C^b_u(G/H_\Rsh)$ is continuous, that is,
$\sigma(G/H_\Rsh)$ is a topological $G$-space. 
The uniformly continuous
mapping of compactification $G/H_\Rsh\to\sigma(G/H_\Rsh)$ has
everywhere dense image and is $G$-equivariant.

It only remains to prove the maximality of $\sigma(G/H_\Rsh)$
as a $G$-equivariant compactification of $G/H$.
Let $X$ be a compact $G$-space, and
let $\phi\colon G/H\to X$ be a continuous $G$-equivariant mapping.
It determines a morphism of $C^\ast$-algebras, $\phi^\ast$, from
$C(X)$ to $C^b_u(G/H_\Rsh)$ via 
\begin{equation}
C(X)\ni f\mapsto [(xH)\mapsto \tilde f(xH):= f(\phi(xH))]\in 
C^b_u(G/H_\Rsh).
\end{equation}
The dual continuous mapping $f^\sim\colon \sigma(G/H_\Rsh)\to
X$ between the Gelfand spaces of the corresponding
$C^\ast$-algebras 
is $G$-equivariant and its restriction to $G/H$ is easily
seen to coincide with $f$. The proof is thus finished.
\end{proof}

\begin{corol} The pair
$(G/H_\Rsh, G)$ has the concentration property if and only
if $\alpha_G(G/H)$ has a fixed point.
\qed
\end{corol}

The superscripts {\it `u'} and {\it `s'} will denote the
uniform (respectively strong) operator topology on the unitary group.
Since the sphere $\s_{\mathcal H}$ is both uniformly and as a
$U({\mathcal H})$-space isomorphic to 
$(U({\mathcal H})_u/\operatorname{St}_\xi)_\Rsh$, 
where $\xi\in\s_{\mathcal H}$ is any and $\operatorname{St}_\xi$
is the stabilizer of $\xi$, we obtain:

\begin{corol} The maximal $U({\mathcal H})_u$-compactification
of the unit sphere of a Hilbert space $\mathcal H$ has 
no fixed points. \qed
\label{fp}
\end{corol}

\begin{remark} One should compare this result with
Stoyanov's theorem \cite{S}: the maximal 
$U({l_2})_s$-compactification of $\s^\infty$
coincides with the unit ball of $l_2$ with the 
weak topology, and thus has a fixed point. Another way to reformulate
Stoyanov's result is of course this: the homogeneous space 
$(U({\mathcal H})_s/\operatorname{St}_\xi)_\Rsh$ is uniformly
isomorphic to the sphere $\s_\H$ equipped with the restriction
of the additive uniform structure of $\H_w$, where the latter denotes
the Hilbert space $\H$ with its {\it weak topology.} 
\tri
\end{remark}

A topological group $G$ is called
{\it extremely amenable (e.a.)} \cite{Pa,P1,P2} if every continuous
action of $G$ on a compact space has a fixed point. 
This property is equivalent to the existence of a fixed point in the
{\it greatest ambit} ${\mathcal S}(G)$ of $G$
\cite{Aus,dV,Te,P2}, that is, the Samuel
compactification of $G_\Rsh$.

\begin{corol} A topological group $G$ is extremely amenable if and
only if the pair $(G_\Rsh, G)$ has the concentration property.
\qed
\end{corol}

\begin{remark}
There are very few known examples of extremely amenable 
topological groups: those due to Herer and Christensen
\cite{HC}, Gromov and Milman \cite{GrM},
Glasner \cite{Gl} (who notes that examples of the same
kind were independently discovered by Furstenberg and B. Weiss but
never published), and the present author \cite{P1}. Nevertheless,
they include
some very natural topological groups having importance in Analysis,
for example the group $U(\infty)$ equipped with the Hilbert-Schmidt 
metric \cite{GrM} (and therefore the groups 
$U(\infty)_u$ and $U(l_2)_s$), and the group $\operatorname{Homeo}_+(I)$ of
orientation-preserving homeomorphisms of the closed interval \cite{P1}.
\tri
\end{remark}

Recall that the Calkin group is the topological factor-group of
$U(l_2)_u$ by the closure of $U(\infty)$ (in the uniform
topology). This closure,
$\overline{U(\infty)}$, is a normal subgroup of $U(l_2)$,
consisting of all operators of the form ${\mathbb I}+T$, where
$T$ is compact.
Denote by $\mathbb T$ the subgroup of $U(l_2)$ consisting of all
scalar multiples of the identity $\lambda{\mathbb I}$, where
$\lambda\in\C$ and $\vert\lambda\vert=1$. 
The group ${\mathbb T}\cdot \overline{U(\infty)}$ forms the
largest proper closed normal subgroup of $U(l_2)_u$ \cite{Ka}.
(Actually, the statement
remains true even without the word `closed'
\cite{dlH2}.) The topological factor-group 
$U(l_2)_u/({\mathbb T}\cdot \overline{U(\infty)})$ is called
the {\it projective Calkin group.}

\begin{corol} 
\label{cal}
The projective Calkin group admits an effective minimal
action on a compact space.
\end{corol}

\begin{proof} By
Proposition \ref{max}, the action of $U({\mathcal H})_u$
upon $\sigma(\s_{\mathcal H})$ is continuous for every Hilbert space
$\mathcal H$, that is, $\sigma(\s_{\mathcal H})$ forms a 
compact $U({\mathcal H})_u$-space. According to a 
result by Gromov and Milman that we cited in the 
Introduction (\cite{GrM}, Example 5.1),
if a compact group $G$ acts by isometries on the unit sphere 
$\s^\infty$ of $l_2$, then the pair
$(\s^\infty,G)$ has the concentration property. It means that there
exists a $\mathbb T$-fixed point $x_1\in\sigma(\s^\infty)$.
Denote by $\mathfrak X$ the closure of the $U(l_2)$-orbit of $x_1$
in $\sigma(\s^\infty)$. It is a compact $U(l_2)_u$-subspace of
$\sigma(\s^\infty)$. Since $\mathbb T$ is the centre of $U(l_2)$,
every point of $\mathfrak X$ is $\mathbb T$-fixed. (In particular, it
follows that $\mathfrak X$ is 
a proper subspace of $\sigma(\s^\infty)$.)

It is a well-known and easy consequence of Zorn's lemma that
every compact $G$-flow contains a minimal subflow (that is,
a non-empty compact $G$-subspace such that the orbit of each point is
everywhere dense in it, see e.g. 
\cite{Aus}). Denote by $\mathcal M$ any minimal subflow of $\frak X$.
Since $U(l_2)$ has no fixed points in 
$\sigma(\s^\infty)$ (Corol. \ref{fp}), 
it follows that every minimal $U(l_2)$-subflow
of $\sigma(\s^\infty)$ is nontrivial, that is, contains more than
one point. In particular, this applies to $\mathcal M$. 

By force of the extreme amenability of the group $U(\infty)_u$ 
(combine the results of \cite{GrM} and \cite{Gl}), 
there is a $U(\infty)$-fixed point, $x^\ast$,
in $\mathcal M$. It follows from the continuity
of the action that $x^\ast$ is also
a fixed point for $\overline{U(\infty)}$, that is,
the stabilizer $\operatorname{St}_{x^\ast}$ contains 
$\overline{U(\infty)}$. Since every point of $\mathcal M$ is 
$\mathbb T$-fixed, it follows that $x^\ast$ is fixed under the action of
the group ${\mathbb T}\cdot\overline{U(\infty)}$.

The stabilizers of
elements of the orbit of $x^\ast$ under the action of
$U(l_2)$ are conjugate to $\operatorname{St}_{x^\ast}$. Since
${\mathbb T}\cdot \overline{U(\infty)}$ is normal 
in $U(l_2)$, every such stabilizer
contains ${\mathbb T}\cdot\overline{U(\infty)}$. Because of minimality of
$\mathcal M$, the $U(l_2)$-orbit of $x^\ast$ is everywhere dense in
$\mathcal M$, and we conclude: all points of $\mathcal M$ are fixed
under the action of ${\mathbb T}\cdot\overline{U(\infty)}$.
It implies that
the action of $U(l_2)_u$ on $\mathcal M$ factors through
an action of the projective Calkin group 
$U({l_2})_u/({\mathbb T}\cdot\overline{U(\infty)})$,
and the latter action is continuous.
Moreover, it is also minimal.

Denote by $K$ the set of all $u\in U(l_2)_u$ leaving each element of
$\mathcal M$ fixed. This is a closed normal subgroup of $U(l_2)_u$, 
containing $\overline{U(\infty)}$, and since it is proper
(in view of minimality and nontriviality of $\mathcal M$), it must be
contained in ${\mathbb T}\cdot \overline{U(\infty)}$ and consequently
coincide with it. It means that the action of the
projective Calkin group is on the compact space $\mathcal M$ is minimal and
effective, and the statement is proved.
\end{proof}

\begin{remark} Contrary to what was in effect claimed in \cite{GrM},
Remark 3.5, the concentration of measure on finite permutation groups 
\cite{Ma} (cf. also \cite{Ta}) does not lead to the
extreme amenability of the infinite
symmetric group $S_\infty$. 
In fact, the group $S_\infty$ of all (finite) permutations of a countably
infinite set $\omega$, equipped with the topology of pointwise 
convergence on the (discrete) set $\omega$, acts effectively on
its universal minimal compact flow and, in particular,
admits continuous actions on compacta without fixed points
(\cite{P1}, Th. 6.5). 
This result, combined with a theorem of Gaughan \cite{Gau}
that every Haudsorff group topology on $S_\infty$ contains the
topology of pointwise convergence,
immediately implies that there is no Hausdorff
group topology making $S_\infty$ into
an extremely amenable topological group.

In particular, $S_\infty$ cannot be made into a
L\'evy group in the sense of \cite{GrM,M1,M2}.
In other words, the concentration of measure 
on the family of finite symmetric groups $S_n$
cannot be observed with respect to a 
right-invariant metric generating a group topology on $S_\infty$.  

The Hamming distance on finite groups of permutations
$S_n$, $n\in\N$ is given by
\begin{equation} 
d(\sigma_1,\sigma_2)=\left\vert\{i\colon \sigma_1(i)\neq
\sigma_2(i)\}\right\vert.
\end{equation}
While the group $S_\infty$ can be represented as the union of an increasing
chain of finite permutation groups $S_n$, the above observation
essentially says that there is no `coherent' way of
putting together the normalised Hamming distances so as to obtain a
right-invariant metric on $S_\infty$. 

In \cite{GrM}, Remark 3.5, it was suggested to consider with
that purpose the function
\begin{equation}
\varphi(\sigma,\eta)=\cases \frac{d(\sigma,\eta)}{\max\{d(\sigma,e),
d(\eta,e)\}}, & \mbox{if $\sigma\neq\eta$,} \\
0 & \mbox{otherwise,}
\endcases
\end{equation}
and then to choose a metric, $\hat d$, on $S_\infty$, determining the
topology of the latter group and Lipschitz equivalent to $\varphi$ with
Lipschitz constant $2$.

Such a metric of course does exist.
However, what matters for the concentration property
and the existence of fixed points in compactifications, 
is not the topology of a topological group $G$ 
{\it per se}, but the uniform structure ${\mathcal U}_\Rsh$
of $G$. Let us show that the uniform structure generated by
$\hat d$ has the property that the right translations of $S_\infty$
do not form a right equicontinuous family,
and therefore this uniform structure does not coincide
with the uniform structure ${\mathcal U}_\Rsh$ 
of {\it any} group topology on $S_\infty$.
(Notice that if $(x,y)\in V_\Rsh$ and $g\in G$, then
$(xg)(yg)^{-1}=xgg^{-1}y^{-1}=xy^{-1}\in V$, that is, $(xg,yg)\in V_\Rsh$
as well, hence the equicontinuity property for right translations.)

In view of the Lipschitz equivalence
of $\hat d$ and $\varphi$, the following sets form a basis of entourages
of the diagonal for
the uniformity generated by $\hat d$
as $\e$ runs over all positive reals:
\begin{equation}
V_\e:=\{(\sigma,\eta)\in S_\infty\colon \varphi(\sigma,\eta)<\e\}.
\end{equation}
Equicontinuity of right translations means that for every
$\e>0$ there exists a $\delta>0$ such that whenever
$(\sigma,\eta)\in V_\delta$ and $\theta\in S_\infty$, one has
$(\sigma\theta,\eta\theta)\in V_\e$. Let $n$ be even, and set
\begin{eqnarray}
\sigma&=&\left(\matrix 1 & 2 & 3 & 4 &5 & 6 &\dots & n-1 & n \\ 
                       2 & 1 & 4 & 3 &6&5&\dots & n & n-1\endmatrix\right), 
                       \nonumber \\
\eta&=&\left(\matrix 1 & 2 & 3 & 4 & 5 & 6&\dots & n-1 & n \\ 
                     1 & 2 & 4 & 3 & 6 & 5 &\dots & n & n-1 \endmatrix\right).
\end{eqnarray}
One has $\varphi(\sigma,\eta)=2/n$ and thus, by choosing $n$ sufficiently
large, we can make the pair $(\sigma,\eta)$ belong to any entourage
$V_\delta$, $\delta>0$. At the same time,
$\varphi(\sigma\eta,\eta^2)=\varphi\left(\left(\matrix
1&2 \\ 2 & 1
\endmatrix\right),e \right)= 2/2=1$, that is, 
the right translation of every entourage of the form
$V_\delta$ is not a subset of $V_{1}$, however small $\delta>0$ be. 
 \tri
\label{sym}
\end{remark}

\section{\label{some}Some lemmas on the geometry of spheres}

Recall that a {\it probabilistic metric space} is a triple,
$(X,\rho,\mu)$, formed by a metric space $(X,\rho)$ and a
normalised ($\mu(X)=1$) Borel measure on $X$.

For a subset $A\subseteq X$ we denote by ${\mathcal O}_\e(A)$
the $\e$-neighbourhood of $A$ in $X$.

The {\it concentration function,}
$\alpha=\alpha_X$, of a probabilistic metric space $X$ is defined 
for each $\e>0$ by
\begin{equation} 
\alpha_X(\e)=1-\inf\left\{\mu\left({\mathcal O}_\e(A)\right)
\colon A\subseteq X \mbox{ is Borel and }
\mu(A)\geq \frac 1 2\right\}
\end{equation}
and $\alpha_X(0)=1/2$. It is a decreasing function in $\e$.

The following observations are straightforward yet useful.
(See e.g. Lemma 3.2 in \cite{M1}.)

\begin{lemma} 
Let $X$ be a probabilistic metric space with the
concentration function $\alpha$.
Let $A\subseteq X$ and $\e>0$. 
\begin{enumerate}
\item If
$\mu(A)>\alpha(\e)$, then $\mu({\mathcal O}_\e(A))> 1/2$.
\item
If $\mu(A)>\alpha(\e/2)$, then $\mu({\mathcal O}_\e(A))\geq 1-\alpha(\e/2)$.
\qed
\end{enumerate}
\label{half}
\end{lemma}

A family $(X_n)_{n=1}^\infty$ of probabilistic metric spaces is called
a {\it L\'evy family} if for each $\e>0$, $\alpha_{X_n}(\e)\to 0$ as
$n\to\infty$,
and a {\it normal L\'evy family} (with constants $C_1,C_2>0$) if for
all $n$ and $\e>0$
\[\alpha_{X_n}(\e)\leq C_1e^{-C_2\e^2n}.\]

By $\mu_n$ we will denote the (unique)
normalized rotation-invariant Borel measure
on the $n$-dimensional Euclidean sphere $\s^n$. The distances 
between points on the spheres will be geodesic distances. In such a way,
the sphere $\s^n$ becomes a probabilistic metric space.
The family of spheres ${\Bbb S}^{n+1}$, $n\in\N$ is normal L\'evy
with constants $C_1=\sqrt{\pi/8}$ and $C_2=1/2$.

For the major concepts, examples, and results
of the theory of concentration of measure on
high-dimensional structures, see 
\cite{FLM,GrM,Ma,M1,M2,MS,St,Ta}.

\begin{lemma} 
There are absolute constants $C_1,C_2>0$ with the following property.
Let $\H$ be a real Hilbert space, let $n\in\N$ and $\e>0$.
Let $P_1$ and $P_2$ be two rank $n$ orthogonal projections in
$\H$, satisfying 
\begin{equation}
\norm{P_1-P_2}_1<n\e,
\end{equation}
where
$\norm\cdot_1$ denotes the trace class operator norm.
Denote by $\s_i$ the unit sphere in the space
$P_i(\H)$, $i=1,2$. Then 
\begin{equation}
\mu_n\left(\{x\in\s_1\colon \norm{x-P_2x}<\e\}\right)
\geq 1-C_1\exp(-C_2\e^2 n).
\end{equation}
\label{const}
\end{lemma}

\begin{proof} Let $0<\delta<1$ be arbitrary and fixed.
Let for every natural
number $n$, $P_1^{(n)}$ and $P_2^{(n)}$ be two rank $n$
projections in $\H$. Assume for a while that, as $n\to\infty$,
\begin{equation}
\mu_n\left(\{x\in\s_1\colon \norm{x-P_2x}<\delta\}\right)
=O(1) \exp(-O(1)n),
\end{equation}
that is, for some positive constants $C_1,C_2>0$ and all $n$ 
one has
\begin{equation}
\mu_n\left(\{x\in\s_1\colon \norm{x-P_2x}<\delta\}\right)\leq
C_1\exp(-C_2n),
\end{equation}
where $\s_i^{(n)}$ denotes the unit sphere in the space 
$P_i^{(n)}(\H)$, $i=1,2$.

Then a standard argument (cf. e.g. \cite{M2}, p. 276) implies that, for
$n$ sufficiently large, there is an orthonormal basis in $P_1^{(n)}(\H)$
formed by elements $e_1,e_2,\dots,e_n\in\s_1^{(n)}$ each
at a distance $>\delta$ from $P_2^{(n)}(\H)$. 
Indeed, fix a $\xi\in\s_1^{(n)}$ and denote by
$\nu$ the Haar probability measure on $SO(n)$. Then
\begin{eqnarray}
\nu\left\{u\in SO(n)\colon \norm{u\xi-P_2(u\xi)}\geq\delta\right\}&=&
\mu_n\left(\{x\in\s_1^{(n)}\colon\norm{x-P_2x}\geq\delta\}\right)
\nonumber \\
&\geq& 1-C_1\exp(-C_2n),
\end{eqnarray}
therefore for any finite subset $F\subseteq\s_1^{(n)}$
\begin{equation}
\nu\left\{u\in SO(n)\colon \norm{x-P_2x}\geq\delta ~~\mbox{for all $x\in uF$}
\right\}
\geq 1-\vert F\vert\cdot C_1\exp(-C_2n).
\end{equation}
If the size $\vert F\vert$ of $F$ grows slower than $C_1\exp(C_2 n)$,
for example if $\vert F\vert$ is polynomial in $n$, then for $n$
satisfying the condition $n<C_1\exp(C_2n)$ we can find
a rotation $u$ taking $F$ outside of the $\delta$-neighbourhood of 
$P_2^{(n)}(\H)$
in $\s_1^{(n)}$. Applying this observation to an arbitrary orthonormal basis
of $P_1^{(n)}(\H)$ as $F$, we obtain an orthonormal system 
$e_1,e_2,\dots,e_n\in\s_1^{(n)}$ with the desired property.

Now extend the collection of $e_i$, $i=1,2,\dots,n$
to an orthonormal basis $(e_i)_{i<\alpha}$
of $\H$. One has
\begin{eqnarray}
\norm{P_1^{(n)}-P_2^{(n)}}_1&=&
\operatorname{Tr}\left(\abs{ P_1^{(n)}-P_2^{(n)} }\right)
\nonumber \\
&\geq& 
\sum_{i<\alpha}\left\vert\left((P_1^{(n)}-P_2^{(n)})(e_i),e_i\right)\right\vert
\nonumber \\
&\geq&
\sum_{i=1}^n\left\vert 1-(P_2^{(n)}(e_i),e_i)\right\vert
\nonumber \\
&>&n\left(1-\sqrt{1-\delta^2}\right).
\label{latter}
\end{eqnarray}
If now $\e>0$ and $\delta>\sqrt{2\e}$, then $1-\delta^2<(1-\e)^2$ and
the latter expression in
(\ref{latter}) is $>n\e$.

The above argument establishes the following:
for each pair of positive constants
$C_1,C_2>0$, if $n$ is so large that $n<C_1\exp(C_2 n)$ and
$P_1,P_2$ are two rank $n$ projections in $\H$ satisfying
$\norm{P_1-P_2}_1<n\e$, then the measure of the set of
points $x\in\s_1$ that are at a distance $<\sqrt{2\e}$ from
$P_2^{(n)}(\H)$ is greater than $C_1\exp(-C_2 n)$.

Since
\begin{equation}
\alpha_{\s^{n+1}}(\e)\leq 
\sqrt{\pi/8}\exp(-\e^2n/2)
\end{equation}
and therefore
\begin{equation}
\alpha_{\s^{n+1}}(\sqrt{2\e})\leq 
\sqrt{\pi/8}\exp(-\e n),
\end{equation}
one has in particular $\mu_n\left(\{x\in\s_1\colon \norm{x-P_2x}
<\sqrt{2\e}\}\right)> \alpha_{\s^n}(\sqrt{2\e})$ whenever $n$ is so
large that 
\begin{equation}
n< \sqrt{\pi/8}\exp(-\e (n-1)).
\label{cond}
\end{equation}

Now observe that the set
\[{\mathcal O}_2:=\{x\in\s_1\colon \norm{x-P_2x}<2\sqrt{2\e}\}\]
contains the open
$\sqrt{2\e}$-neighbourhood of the set 
\[{\mathcal O}_1:=\{x\in\s_1\colon \norm{x-P_2x}<\sqrt{2\e}\}.\]

According to Lemma \ref{half}, (2), one has
\begin{eqnarray}
\mu_n\left(\{x\in\s_1\colon \norm{x-P_2x}
<\sqrt{2\e}\}\right)&\geq& 1-\alpha_{\s^n}(\sqrt{2\e}) \nonumber \\
&\geq&
1-\alpha_{\s^{n+1}}(\sqrt{2\e}) \nonumber \\
&\geq&
1- \sqrt{\frac{\pi}{8}}e^{-\frac 1 4 \e n}
\end{eqnarray}
whenever $n$ satisfies (\ref{cond}).
Replacing $\e$ in both formulae with $\e^2/2$ yields the following: 
\begin{equation}
\mu_n\left(\{x\in\s_1\colon \norm{x-P_2x}
<\e\}\right)\geq 
1- \sqrt{\frac{\pi}{8}}e^{-\frac 1 8 \e^2 n}
\end{equation}
whenever $n$ is so large that
\begin{equation}
n< \sqrt{\frac\pi 8}\exp\left(-\frac 1 2 \e^2 (n-1)\right)
\end{equation}
and $P_1,P_2$ are two orthogonal projections of rank $n$ satisfying
$\norm{P_1-P_2}_1<n\e$.
And this is the desired result in slight disguise.
\end{proof}

\begin{remark}
Of course, in general one does not expect 
{\it all} points of $\s_1$ to be at a distance $<\e$
from $\H_2$. 
Consider the projections $P_A$, $P_B$ in the space
$l_2$, where $A,B$ are two distinct subsets 
of the index set $\N$ having the same finite
cardinality $n$. If $i\in A\setminus B$,
then $e_i\in \s_1$ and $d(e_i,\H_2)=1$. At the same time,
$P_A$ and $P_B$ can be chosen as to satisfy the condition
$\norm{P_A-P_B}_1<n\e$ with $\e>0$ is as small as desired, as
the `F\o lner ratio' $\vert A\Delta B\vert/\vert A\cup B\vert \to 0$.
\tri
\end{remark}

\begin{lemma} Let $P_1$ and $P_2$ be projections of the same finite
rank $n$ in a (real or complex)
Hilbert space $\H$. Then there exist one-dimensional
projections $e_j^i$, $j=1,2$, $i=1,2,\dots,n$, such that
\begin{equation}
P_j=\vee_{i=1}^n e_j^i,~~ j=1,2
\end{equation}
and
\begin{equation}
e_j^i\perp e_k^m ~~\mbox{ whenever $j,k\in\{1,2\}$,
$i,m\in\{1,2,\dots, n\}$, and $i\neq m$.}
\end{equation}
\label{projections}
\end{lemma}

\begin{proof}
Let $x$ be an eigenvector of $P_1+P_2$ corresponding to an
eigenvalue $\lambda$. 
Then $\lambda= 1\pm \cos\theta$, where
$\theta\in [0,\pi/2]$ is the angle between one-dimensional subspaces
spanned by $P_1x$ and $P_2x$. The space $(P_1\vee P_2)(\H)$
has an orthogonal basis formed by 
eigenvectors of $P_1+ P_2$ which can be written in the form
$x_i^\pm$, $i=1,2,\dots,n$, where $x_i^\pm$ corresponds to the 
eigenvalue $\lambda_i^\pm=1\pm\cos\theta_i$, $\theta_i$ being as
above; if $\theta_i=0$, then we only consider $\lambda^+_i=2$
and $x_i^+$.
Since $P_j(x^+_i)$, $j=1,2$ span the same subspace of dimension
$2$ (respectively $1$ where $\theta_i=0$)
as the vectors $x_i^+$ and $x_i^-$ do (respectively, the vector
$x_i^+$ does), it follows that
$P_j(x^+_i)\perp P_k(x^+_m)$ whenever $i\neq m$.
Now let $e_j^i$ be the projection onto the one-dimensional
subspace spanned by $P_j(x^+_i)$.
\end{proof}

\begin{lemma} There are absolute constants $C_1',C_2'>0$ with the
property that,
under the assumptions and using the notation of Lemma
\ref{const}, there exists an isometry $r\colon\s_1\to\s_2$
with
\begin{equation}
\mu_n\left(\{x\in\s_1\colon \norm{x-r(x)}<\e\}\right)
\geq 1-C'_1\exp(-C'_2\e^2 n).
\end{equation}
\label{iso}
\end{lemma}

\begin{proof} Let the projections $e_j^i$, $j=1,2$, $i=1,2,\dots,n$
be as in Lemma \ref{projections}. For every $i=1,2,\dots,n$ denote by
$r_i$ the (unique) isometric isomorphism from the 
one-dimensional range of $e_1^i$ to that of $e_2^i$.
(That is, the reflection across the linear span of $x_i^+$, using the
notation from the proof of Lemma \ref{projections}.)
Since $P_2P_1x_i^+=(\lambda-1)P_2x_i^+\in\operatorname{span}(P_2x_i^+)$, 
one concludes that
if $x\in \operatorname{span}(P_1x_i^+)$ 
and $\norm x=1$, then
\begin{equation}
\norm{r_i (x)-x}\leq \sqrt 2\operatorname{dist}(x,P_2(x_i^+))=
\sqrt 2\operatorname{dist}(x,\H_2).
\label{estimates}
\end{equation}
The orthogonal sum of linear operators
$r=\oplus_{i=1}^n r_i$ is an isometry between
$\H_1$ and $\H_2$. The equation (\ref{estimates}) implies
that for each $x\in\H_1$, $\norm x=1$, one has
\begin{equation}
\norm{r (x)-x}\leq \sqrt 2\operatorname{dist}(x,\H_2)
\leq\sqrt 2\operatorname{dist}(x,\s_2),
\end{equation}
and the proof is finished by applying Lemma \ref{const}.
\end{proof}

\begin{remark}
\label{complex}
As an immediate corollary of the statement in the real case,
Lemma \ref{iso} remains true in a complex Hilbert space $\H$
as well.
\tri
\end{remark}

\section{\label{fromam}From amenability to concentration property}

From the previous work \cite{GrM,M1,M2,P3} we know that if a group
$G$ is compact or discrete amenable, and $\pi$ is
a unitary representation of $G$ in an infinite-dimensional Hilbert
space, then the dynamical system $(\s_\pi,G)$ has the property
of concentration. Our present aim is to push this result further
as far as possible. A plausible-looking conjecture might be that
the conclusion remains in force if $\pi$ is just an amenable representation
of a group in an infinite-dimensional Hilbert space. 
However, this is not true.

\begin{example} 
\label{ex2}
Let $G=F_2$ be the free group on two generators. Denote by
$\pi_1$ an irreducible unitary representation of $F_2$ in a
finite-dimensional Hilbert space $\H_1$, and let
$\pi_2$ denote the left regular representation
of $F_2$ in the space $\H_2=l_2(F_2)$. Let $\pi=\pi_1\oplus\pi_2$
be the direct sum representation of $F_2$ in the Hilbert space
$\H=\H_1\oplus\H_2$.
Since $\pi$ contains a finite-dimensional subrepresentation
$\pi_1$, it is amenable (\cite{B}, Th. 1.3, (i)+(ii)). Both $\pi_1$ and
$\pi_2$ do not have the concentration property (cf. Ex. \ref{expl}).
Now it is enough to apply Proposition \ref{sum}.
\tri
\end{example}

Nevertheless, the above situation --- in which amenability of a
representation only stems from the presence of a single 
finite-dimensional subrepresentation ---
is, in fact, the only one where the concentration property is
not to be found.
If a representation is amenable in a `nontrivial way,' 
the concentration property of spheres rebounds.
The rest of this section will be devoted to establishing 
the corresponding
result (Theorem \ref{conc}), which not only generalizes all the
previously obtained results in this direction, but is, in a sense,
`at the end of the road.'

\begin{lemma}
Let $F$ be a finite collection of unitary operators on a 
(real or complex) Hilbert
space $\H$. Suppose that for every $\e>0$ and every natural 
$k$ there is an orthogonal projection
$P$ in $\H$ of rank $n\geq k$ such that
\begin{equation}
\norm{gPg^{-1}-P}_1<n\e
\end{equation}
for all $g\in F$. Then the system $(\s_\H,F)$ has
the concentration property.
\label{pro}
\end{lemma}

\begin{proof} 
Choose a sequence of orthogonal projections
$P_n$, $n\in\N$, having the properties:
\begin{enumerate}
\item
$r_n=\operatorname{rank} P_n\to\infty$, and
\item
$\norm{gP_ng^{-1}-P_n}_1<r_n/n$ for all $g\in F$ as $n\to\infty$.
\end{enumerate}
Denote by $\s_n=\s^{r_n}$ the unit sphere in the space
$P_n(\H)$, and for each $g\in F$ denote by $\s_n^g=g\s_n$
the unit sphere in the $r_n$-dimensional
space $gP_n(\H)\equiv gP_ng^{-1}(\H)$.

Now let $\gamma=\{A_1,A_2,\dots,A_m\}$ be an arbitrary finite cover
of the unit sphere $\s_\H$. 
Clearly, for at least one $i=1,2,\dots,m$
the set of natural numbers
\begin{equation}
\{n\in\N\colon \mu_{\s_n}(A_i\cap \s_n)\geq 1/m\}
\end{equation}
is infinite. 
We claim that the set $A=A_i$ is then $F$-essential.

Proceeding to a subsequence of the selected sequence of projections if
necessary, we can assume without loss in generality that
$\mu_{\s_n}(A_i\cap \s_n)\geq 1/m$ for all $n\in\N$.
For every $g\in F$ the measure of the set $gA\cap \s^g_n$ in
the latter sphere is the same as the measure of $A\cap \s_n$,
and therefore $\geq 1/m$ for all $n$. 

Let an $\e>0$ be arbitrary.
According to L\'evy's concentration of
measure property of spheres, the measure of every set of the form
$g{\mathcal O}_\e(A)\cap\s^g_n$ in $\s^g_n$, where 
$g\in F\cup\{e\}$, is $1-O(1)\exp(-O(1)\e^2n)$. 

An application of Lemma \ref{iso} and Remark \ref{complex} 
to $P=P_1$ and 
$gPg^{-1}=P_2$ yields that for every $g$ as above the
measure of the set $g{\mathcal O}_{2\e}(A)\cap\s_n$ in $\s_n$
is $1-O(1)\exp(-O(1)\e^2n)$. Indeed, there is an isometry
$i_g\colon\s^g_n\to\s_n$ with the property that for all points $x$
of $\s^g_n$ apart from a set of measure $O(1)\exp(-O(1)\e^2n)$
one has $\norm{x-i_g(x)}<\e$. Consequently, the set
$g{\mathcal O}_{2\e}(A)\cap\s_n$ contains the set 
$i_g(g{\mathcal O}_\e(A)\cap\s^g_n)$. The measure of the latter set
in $\s_n$ is $1-O(1)\exp(-O(1)\e^2n)$, because so is the measure of
$g{\mathcal O}_\e(A)\cap\s^g_n$ in $\s^g_n$ and
the isometry
$i_g$ between the spheres is automatically a measure-preservig map.

We conclude that if $n$ is sufficiently large, then the
sets $g{\mathcal O}_{2\e}(A)$, $g\in F$, 
have a non-empty common part, and indeed the measure
of its intersection with the sphere $\s_n$ is
$1-O(1)\abs{F}\exp(-O(1)\e^2n)$. This finishes the proof.
\end{proof}

\begin{lemma} 
\label{rank} Let $\pi$ be a unitary representation of a group $G$ in
a Hilbert space $\H$. If $\pi$ is not of the form
$\pi_1\oplus\pi_2$, where $\pi_1$ is finite-dimensional and
$\pi_2$ is non-amenable, then for every finite subset $F\subseteq G$,
every $\e>0$ and every natural $k$ there is an orthogonal projection
$P$ in $\H$ of rank $n\geq k$ such that
\begin{equation}
\norm{\pi(g)P\pi(g)^{-1}-P}_1<n\e
\end{equation}
for all $g\in F$. 
\label{high}
\end{lemma}

\begin{proof}
If $\pi$ admits finite-dimensional subrepresentations of
arbitrarily high dimension, then the desired projections can be
constructed in an obvious way. Otherwise, using
the assumption of the lemma, one can
assume without loss in generality till the end of
the proof that $\pi$ contains no nontrivial
finite-dimensional subrepresentations.

According to the F\o lner property of amenable representations
as established by Bekka
(\cite{B}, Th. 6.2), for every finite $F\subseteq G$ and
each $\e>0$ there is an orthogonal projection
$P=P_{F,\e}$ in $\H$ of finite rank such that
\begin{equation}
\norm{\pi(g)P\pi(g)^{-1}-P}_1<\norm P_1\e
\label{property}
\end{equation}
for all $g\in F$. 

Suppose it is not in general possible to choose such a $P$ of an
arbitrarily high finite rank. In such a case, there are a finite
$\Phi\subseteq G$ and 
an $\e'>0$ such that for each finite $F\supseteq \Phi$ and
each $\e<\e'$ an arbitrary projection 
$P$ satisfying (\ref{property})
has rank $\norm P_1\leq N$, where $N$ is a fixed natural number.
Moreover, one can assume without loss in generality that the
equality is always achieved for a suitable $P=P_{F,\e}$.

Notice that on the collection of all projections of fixed finite rank
$N$ the trace class metric and the operator metric are both Lipschitz
equivalent to the Hausdorff distance between unit spheres
in the range spaces of the projections.
(The equivalence of the trace class and operator metrics
follows from the obvious inequality $\operatorname{Tr}(\vert P_1-P_2\vert)
\leq 2N\norm{P_1-P_2}$ for every two projections $P_1,P_2$ of
rank $N$.
For the equivalence of the operator and Hausdorff metrics,
see \cite{O}, 3.12 and 3.4.(h).)
The space of all projections of rank $N$ is thus uniformly
isomorphic to the Grassmannian $\operatorname{Gr}_N(\H)$ and forms a
complete uniform space. (Cf. \cite{O}, 3.12.(c) and 3.4.(e).)

For every $F$ and $\e$ as above denote by ${\mathcal P}_{F,\e}$ the
non-empty set of all projections of rank exactly
$N$ satisfying (\ref{property}).
Now equip the set $\mathcal{P}_\omega(G)\times\R_+$
of all pairs $(F,\e)$ as above with the product
partial order making it into a directed set. 
The diameters of ${\mathcal P}_{F,\e}$ cannot converge to zero
over $\mathcal{P}_\omega(G)\times\R_+$. 
Otherwise the sets  ${\mathcal P}_{F,\e}$
would form a Cauchy prefilter having
a limit point $P$, which is again a projection of 
rank $N$ satisfying the property (\ref{property}) for every
finite $F\subseteq G$ and each
$\e>0$. In other words, $P$ commutes with every $g\in G$, that is,
$P(\H)$ is the space of a nontrivial finite-dimensional subrepresentation
of $\pi$, leading to a contradiction.
We conclude that for some $\delta>0$, the set of pairs $(F,\e)$
satisfying
\begin{equation}
\operatorname{diam}\left({\mathcal P}_{F,\e} \right)\geq\delta
\end{equation}
is cofinal in $\mathcal{P}_\omega(G)\times\R_+$.

It remains to notice that, if an arbitrary finite set $F\subseteq G$ 
is fixed, then having at one's disposal,
for every $\e>0$, a pair of
projections $P_1$ and $P_2$ of the same finite rank $N$ 
satisfying 
\begin{equation}
\norm{\pi(g)P_i\pi(g)^{-1}-P_i}<\e
\label{pred}
\end{equation}
for all $g\in F$ and $i=1,2$,
and also the condition
\begin{equation}
\norm{P_1-P_2}\geq\delta
\end{equation}
for a fixed $\delta>0$ enables one to produce a new projection
$P'$ of rank $\geq N+1$ satisfying 
\begin{equation}
\norm{\pi(g)P'\pi(g)^{-1}-P'}<\e_1,
\end{equation}
where $\e_1\to 0$ as $\e\to 0$, thus obtaining a contradiction with
the presumed maximality of $N$.

Indeed, choose two sequences of projections $P_i^{(1)}$ and
$P_i^{(2)}$ of rank $N$ having the properties
$\norm{P_i^{(1)}-P_i^{(2)}}\geq\delta$ and
$\norm{\pi(g)P_i^{(j)}\pi(g)^{-1}-P_i}\to 0$ as $i\to\infty$
for all $g\in F$ and $j=1,2$. Using Lemma \ref{projections},
choose for every $i\in\N$ one-dimensional
projections $_ke_j^i$, $j=1,2$, $k=1,2,\dots,N$, such that
\begin{equation}
P_i^{(j)}=\vee_{k=1}^N {\,}_ke_j^i,~~ j=1,2
\end{equation}
and
\begin{equation}
_ke_j^i\perp \,_ke_l^m ~~\mbox{ whenever $j,l\in\{1,2\}$,
$k,m\in\{1,2,\dots, N\}$, and $k\neq m$.}
\end{equation}
Without loss of generality and proceeding to a subsequence if
necessary, one can assume that $\norm{_ke^i_1-{_ke^i_2}}\to c_k$
for each $k$ as $i\to\infty$, where 
$0\leq c_1\leq c_2\leq\cdots\leq c_N$. Let $d$ be the smallest integer
$\leq N$
with the property $c_d>0$; clearly, such a $d$ exists.
Let \[P_i=\vee_{k=1}^{d-1}{_ke_1^i}\vee\vee_{k=d}^N{_ke_1^i}\vee
\vee_{k=d}^N{_ke_2^i}.\]
The rank of the projection $P_i$ is $2N-d+1>N$.
The space $P_i(\H)$ is the direct sum of subspaces
$\H_i$, $i=1,2,\cdots,N$, where
$\H_i=\operatorname{span}({_ke_1^i})$ for $k<d$ and
$\H_i=\operatorname{span}({_ke_1^i},{_ke_2^i})$ for $k\geq d$.

For every $k\geq d$ and every $g\in F$, 
the Hausdorff distance between the unit spheres
$\s_{\H_k}$ and $\pi(g)\s_{\H_k}\pi(g)^{-1}$ approaches zero
as $i\to\infty$, because (i) the geodesic distances between 
${_ke_1^i}$ and ${_ke_2^i}$ are bounded from below 
by some positive constant and at the same time do not exceed $\pi/2$,
and (ii) the distances between ${_ke_j^i}$ and
$\pi(g){_ke_j^i}\pi(g)^{-1}$, $i=1,2,\dots,N$, $j=1,2$ converge to zero
as $i\to\infty$.
For $k<d$ the same is true in a trivial sort of way.
Since the spheres $\s_{\H_k}$, $k=1,2,\cdots,n$ are paiwise
orthogonal and the same is true of the spheres
$\pi(g)\s_{\H_k}\pi(g)^{-1}$, one concludes that the Hausdorff distance
between the unit sphere $\s_i$ in the space $P_i(\H)$ and
the unit sphere $\pi(g)\s_i\pi(g)^{-1}$ approaches zero as $i\to\infty$.
Since the rank of $P_i$ is bounded above by $2N$,
this means $\norm{P_i-gP_ig^{-1}}_1\to 0$, as required.
\end{proof}

\begin{thm} 
\label{conc} Let $\pi$ be a unitary representation of a group $G$ in
an infinite-dim\-en\-sion\-al Hilbert space $\H$. If every
subrepresentation of $\pi$ having finite codimension is amenable
(that is, $\pi$ is not of the form
$\pi_1\oplus\pi_2$, where $\pi_1$ is finite-dimensional and
$\pi_2$ is non-amenable), then 
the dynamical system $(\s_\H,G,\pi)$ has the concentration property.
\end{thm}

\begin{proof} 
Combine Lemmas \ref{high} and \ref{pro}.
\end{proof}

\section{\label{vis}Amenability {\it vis-\`a-vis} concentration property}
Now we can deduce all our main results.
To begin with, we obtain a description of subgroups 
of the full unitary
group $U(\H)$ whose action on the unit sphere has the concentration
property.

\begin{thm} Let $\pi$ be a unitary representation of a group
$G$ in a Hilbert space $\H$. The system
$(\s_\H,G,\pi)$ has the concentration property if and only if
\begin{itemize}
\item either $\pi$ has a non-zero invariant vector, or
\item $\dim\H=\infty$ and
every subrepresentation of $\pi$ of finite codimension 
is amen\-able.
\end{itemize}
\label{either}
\end{thm}

\begin{proof} $\Rightarrow$: 
If $\dim\H<\infty$, then under our assumption $\pi$ clearly
has a non-zero invariant vector. Otherwise,
suppose there exists a non-amenable 
subrepresentation $\pi_1$ of $\pi$ having finite codimension.
If the finite-dimensional subrepresentation $\pi_1^\perp$ has no
non-zero
invariant vectors, then it does not have the concentration property
and the same is true of $\pi$ according to
Proposition \ref{sum}. $\Leftarrow$: Immediate from Theorem \ref{conc}.
\end{proof}

The following particular cases are of some interest.

\begin{corol}
If a unitary representation of a group $G$ in a Hilbert space $\H$
is amenable and
has no finite-dimensional subrepresentations, then the system
$(\s_\H,G,\pi)$ has the concentration property.
\label{no}
\qed
\end{corol}

\begin{corol} If $\pi$ is an amenable 
irreducible unitary representation of
a group $G$ in an infinite-dimensional Hilbert space $\H$, then 
the system $(\s_\H,G,\pi)$ has the concentration property.
\qed
\end{corol}

Now we are able to extend our criterion of amenability stated in \cite{P3}
from discrete groups to locally compact ones.

\begin{thm}
\label{charact}
A locally compact group $G$ is amenable if and only if for every
strongly continuous unitary representation $\pi$ of $G$ in an
infinite-dimensional Hilbert space $\H$, the dynamical system
$(\s_\H,G,\pi)$ has the concentration property.
\end{thm}

\begin{proof} $\Rightarrow$: every strongly continuous unitary
representation of an amenable locally compact group is amenable
\cite{B}, and therefore so are all its subrepresentations,
and Theorem \ref{either} applies.
$\Leftarrow$: if $G$ is finite, there is nothing to prove, otherwise
$L_2(G)$ is infinite-dimensional and Corollary \ref{l2} applies
together with our assumption.
\end{proof}

Conversely, we can characterize amenable representations in terms of
the concentration property.

\begin{thm}
\label{amenability}
A unitary representation $\pi$ of a group $G$ in a Hilbert space $\H$
is amenable if and only if 
\begin{itemize}
\item 
either $\pi$ contains a finite-dimensional
subrepresentation, or 
\item the $G$-space $(\s_\H,G,\pi)$ has the
concentration property.
\end{itemize}
\end{thm}

\begin{proof} 
$\Rightarrow$: If $\pi$ has no finite-dimensional subrepresentations,
then Corollary \ref{no} applies.
$\Leftarrow$: Combine Corollary \ref{if} with the
following: a unitary representation is amenable if it
contains a finite-dimensional subrepresentation (\cite{B}, Thm. 1.3).
\end{proof}

As an application of our techniques, we show that amenability of
a representation of a group $G$ is equivalent to the existence of a 
L\'evy-type $G$-invariant 
integral for functions on the sphere in the
space of representation. Namely, we are able to invert 
Proposition \ref{inv} and obtain a new equivalent definition of
an amenable representation very much in the classical spirit
of amenability.

\begin{thm}
Let $\pi$ be a unitary representation of a group
$G$ in a Hilbert space $\H$. The following conditions are
equivalent.
\begin{itemize}
\item The space $C^b_u(\s_\H)$ of all 
bounded uniformly
continuous functions on the unit sphere $\s_\H$  
admits a $G$-invariant mean.
\item The representation $\pi$ is amenable.
\end{itemize}
\label{levy}
\end{thm}

\begin{proof} $\Rightarrow$: Prop. \ref{inv}.
$\Leftarrow$: If $\pi$ contains a finite-dimensional subrepresentation
$\pi_1=\pi\vert_{\H_1}$, then the desired $G$-invariant mean 
is obtained by integrating the (restriction of) an $f\in C^b_u(\s_\H)$ 
over the unit sphere of $\H_1$.
If $\pi$ contains no finite-dimensional subrepresentations, then,
according to Corollary \ref{no} and Proposition \ref{equiv}.(\ref{four}), 
there even exists a multiplicative
$G$-invariant mean on $C^b_u(\s_\H)$.
\end{proof}

\subsection*{Acknowledgement}
Thanks to Pierre de la Harpe for his stimulating remarks on my e-print 
\cite{P3}.

\bibliographystyle{amsplain}

\end{document}